\documentclass[12pt]{article}
\usepackage{amsmath}
\usepackage{amsbsy}
\usepackage{amsfonts}
\usepackage{mathtools}
\usepackage{bbm}
\usepackage{amsmath}
\usepackage{amstext}
\usepackage{amsthm}
\usepackage{bbold}
\usepackage{amssymb}
\usepackage[dvips]{graphicx}
\usepackage{amsfonts}
\usepackage{caption2}
\usepackage{amsfonts,mathrsfs}
\usepackage{makeidx}
\usepackage{fancyhdr}
\usepackage{epic,eepic,epsfig}
\usepackage{mdwlist}
\usepackage{amscd}
\usepackage{amsbsy}
\usepackage{epsfig,graphicx}
\usepackage{epsfig,mathrsfs,wasysym,stmaryrd}
\usepackage[all]{xypic}
\usepackage[knot,poly]{xy}
\usepackage{latexsym}
\usepackage{comment}

\theoremstyle{plain}

\newtheorem{theorem}{Theorem}

\newtheorem{lemma}{Lemma}

\newtheorem{corollary}{Corollary}
\newenvironment{pf}{\medskip\noindent{Proof:}
  \hspace{-.5cm}      \enspace}{\hfill \qed \newline \smallskip}

\date{}

\begin{document}

\begin{center}
\textbf{\LARGE{Complexity of circulant graphs with non-fixed
jumps, its  arithmetic properties and asymptotics}}
\vspace{12pt}

{\large\textbf{A.~D.~Mednykh,}}\footnote{{\small\em Sobolev Institute of Mathematics,
Novosibirsk State University, smedn@mail.ru}}
{\large\textbf{I.~A.~Mednykh,}}\footnote{{\small\em Sobolev Institute of Mathematics,
Novosibirsk State University, ilyamednykh@mail.ru}}
\end{center}

\section*{Abstract}

In the present paper, we investigate a family of circulant graphs with non-fixed   
jumps $$G_n=C_{\beta n}(s_1, \ldots,s_k,\alpha_1n,\ldots,\alpha_\ell n),\,
 1\le s_1<\ldots<s_k\le[\frac{\beta n}{2}],\, 1\le \alpha_1< \ldots<\alpha_\ell\le[\frac{\beta}{2}].$$
Here $n$ is an arbitrary  large natural number and integers $s_1, \ldots,s_k,\alpha_1, \ldots,\alpha_\ell$ are supposed to be fixed.  

First, we present an explicit formula for the number of spanning trees in the graph
$G_n.$ This formula
is a product of $\beta s_k-1$ factors, each given by the $n$-th Chebyshev polynomial
of the first kind evaluated at the roots of some prescribed polynomial of degree  $s_k.$
Next, we provide some arithmetic properties of the complexity function. We show that
the number of spanning trees in   $G_n$ can be represented in the form
$\tau(n)=p \,n \,a(n)^2,$ where $a(n)$ is an integer sequence and $p$ is a prescribed
natural number depending of parity of $\beta$ and $n.$
Finally, we find an asymptotic formula for $\tau(n)$ through the Mahler measure of the Laurent polynomials differing by a constant from
$2k-\sum\limits_{i=1}^k(z^{s_i}+z^{-s_i}).$\bigskip

\noindent
\textbf{Key Words:} spanning tree, circulant graph, Laplacian matrix, Chebyshev polynomial,
Mahler measure\\\textbf{AMS classification:} 05C30, 39A10\\

\section{Introduction}

The \textit{complexity} of a finite connected graph $G$, denoted by $\tau(G),$ is the number of spanning trees of $G.$ The famous Kirchhoff's Matrix Tree Theorem~\cite{Kir47} states that $\tau(G)$ can be expressed as the product of non-zero Laplacian eigenvalues of $G$ divided by the number of its vertices. Since then, a lot of papers devoted to the complexity of various classes of graphs were published. In particular, explicit formulae were derived for complete multipartite graphs~\cite{Lew99}, wheels~\cite{BoePro}, fans~\cite{Hil74}, prisms~\cite{BB87}, anti-prisms \cite{SWZ16}, ladders~\cite{Sed69}, M\"obius ladders~\cite{Sed70}, lattices \cite{SW00} and    other families. The complexity of circulant graphs has been the subject of study by many authors \cite{Xiebin, XiebinLinZhang, YTA97, ZY99, ZhangYongGol, ZhangYongGolin, ZhangGolin2002, MZXX2015}.

Starting with Boesch and Prodinger \cite{BoePro} the idea to calculate the complexity of graphs by making use of Chebyshev polynomials was implemented. This idea provided a way to find complexity of circulant graphs and their natural generalisations in~\cite{KwonMedMed, Louis, Med1, Xiebin, ZhangYongGolin, ZhangGolin2002}.

Recently, asymptotical behavior of complexity for some families of graphs was investigated 
from the point of view of so called Malher measure \cite{GutRog, SilWil, SilWil1}. For general properties of the Mahler measure 
see, for example \cite{Smyth08} and \cite{EverWard}. It worth mentioning that the Mahler 
measure is related to the growth of groups, values of some hypergeometric functions and 
volumes of hyperbolic manifolds \cite{Boy02}.

For a sequence of graphs $G_n,$ one can consider the number of vertices $v(G_n)$ and the number of spanning trees $\tau(G_n)$ as functions of $n.$ Assuming that $\lim_{n\to \infty}\frac{\log\tau(G_n)}{v (G_n)}$ exists, it is called the thermodynamic limit of the family $G_n$ \cite{Lyon05}. This number plays an important role in statistical physics and was investigated by many authors (\cite{Wu77}, \cite{SW00}, \cite{Kas2}, \cite{SilWil}, \cite{SilWil1}).

The purpose of this paper is to present new formulas for the number of spanning trees in circulant graphs with non-fixed jumps and investigate their arithmetical properties and asymptotics. We mention that the number of spanning trees for such graphs was found earlier in \cite{GolinXY2010, Louis, XiebinLinZhang, ZhangYongGol, ZhangYongGolin, MZXX2015}. Our results are different from those obtained   in the cited papers. Moreover, by the authors  opinion, the obtained formulas are more convenient for analytical investigation.

The content of the paper is lined up as follows. Basic definitions and preliminary results are given in sections~\ref{basic} and \ref{polynomials}. Then, in the section~\ref{count}, we present a new explicit formula  for the number of spanning trees in the undirected circulant graph 
$$C_{\beta n}(s_1,s_2,\ldots,s_k,\alpha_1n,\alpha_2n,\ldots,\alpha_{\ell} n),\,1\le s_1<\ldots
<s_k\le[\frac{\beta n}{2}],\, 1\le \alpha_1< \ldots<\alpha_{\ell}\le[\frac{\beta}{2}].$$
This formula is a product of $\beta s_k-1$ factors, each given by the $n$-th Chebyshev polynomial of the first kind evaluated at the roots of a prescribed polynomial of degree  $s_k.$ We note the case $\beta=1$ and $\ell=0$  of the circulant graphs with bounded jumps has been investigated in our previous papers \cite{MedMed2, MedMed3}.

Next, in the section~\ref{circarithm}, we provide some arithmetic properties of the complexity function. More precisely, we show that the number of spanning trees of the circulant graph can be represented in the form $\tau(n)=p\,n\,a(n)^2,$ where $a(n)$ is an integer sequence and $p$ is a prescribed natural number depending only of parity of $n.$ Later, in the section~\ref{assection}, we use explicit formulas for the number of spanning trees to produce its asymptotics through the Mahler measures of the finite set of Laurent polynomials 
$P_u(z)=2k-\sum\limits_{i=1}^k(z^{s_i}+z^{-s_i})+4\sum\limits_{m=1}^{\ell}
\sin^2(\frac{\pi\, u \,\alpha_m  }{\beta}),\,u=0,1,\ldots,\beta-1.$ 
As a consequence (Corollary~\ref{mahler}), we prove that the thermodynamic limit of sequence
$C_{\beta n}(s_1,s_2,\ldots,s_k,\alpha_1n,\alpha_2n,\ldots,\alpha_\ell n)$ as $n\to\infty$ is the arithmetic mean of small Mahler measures of Laurent polynomials $P_u(z),\,u=0,1,\ldots,\beta-1.$ In the section~\ref{tables}, we illustrate the obtained results by a series of examples.

\section{Basic definitions and preliminary facts}\label{basic}

Consider a connected finite graph $G,$ allowed to have multiple edges but without
loops. We denote the vertex and edge set of $G$ by $V(G)$ and $E(G),$ respectively.
Given $u, v\in V(G),$ we set $a_{uv}$ to be equal to the number of edges between
vertices $u$ and $v.$ The matrix $A=A(G)=\{a_{uv}\}_{u, v\in V(G)}$ is called
\textit{the adjacency matrix} of the graph $G.$ The degree $d(v)$ of a vertex
$v \in V(G)$ is defined by $d(v)=\sum_{u\in V(G)}a_{uv}.$ Let $D=D(G)$ be the
diagonal matrix indexed by the elements of $V(G)$ with $d_{vv} = d(v).$ The matrix
$L=L(G)=D(G)-A(G)$ is called \textit{the Laplacian matrix}, or simply \textit{Laplacian},
of the graph $G.$

In what follows, by $I_n$ we denote the identity matrix of order $n.$

Let $s_1,s_2,\ldots,s_k$ be integers such that $1\leq s_1, s_2, \ldots, s_k\leq\frac{n}{2}.$
The graph $G=C_{n}(s_1,s_2,\ldots,s_k)$ with $n$ vertices $0,1,2,\ldots,~{n-1}$ is called
\textit{circulant graph} if the vertex $i,\, 0\leq i\leq n-1$ is adjacent to the vertices
$i\pm s_1,i\pm s_2,\ldots,i\pm s_k\ (\textrm{mod}\ n).$ All vertices of
the graph $G$ have even degree $2k.$ If there is $i$  such that $s_i=\frac{n}{2}$ then graph
 $G$  has multiple edges.

We call an $n\times n$ matrix \textit{circulant,} and denote it by $circ(a_0, a_1,\ldots,a_{n-1})$
if it is of the form
$$circ(a_0, a_1,\ldots, a_{n-1})=
\left(\begin{array}{ccccc}
a_0 & a_1 & a_2 & \ldots & a_{n-1} \\
a_{n-1} & a_0 & a_1 & \ldots & a_{n-2} \\
  & \vdots &   & \ddots & \vdots \\
a_1 & a_2 & a_3 & \ldots & a_0\\
\end{array}\right).$$

It easy to see that adjacency and Laplacian matrices for the circulant graph is circulant matrices.
The converse is also true. If the Laplacian matrix of a graph is circulant then the graph is also circulant.

Recall \cite{PJDav} that the eigenvalues of matrix $C=circ(a_0,a_1,\ldots,a_{n-1})$
are given by the following simple formulas $\lambda_j=L(\zeta_n^j),\,j=0,1,\ldots,n-1,$
where $L(x)=a_0+a_1 x+\ldots+a_{n-1}x^{n-1}$ and $\zeta_n$ is an order $n$ primitive
root of the unity. Moreover, the circulant matrix $C=L(T),$ where $T=circ(0,1,0,\ldots,0)$
is the matrix representation of the shift operator
$T:(x_0,x_1,\ldots,x_{n-2},x_{n-1})\rightarrow(x_1, x_2,\ldots,x_{n-1},x_0).$

Let $P(z) = a_0 +a_1 z+\ldots+a_d z^d  = a_d \prod\limits_{k=1}^d(z-\alpha_k)$ be a non-constant polynomial with complex coefficients. Then, following Mahler \cite{Mahl62} its \textit{Mahler measure} is defined to be 
\begin{equation} M(P):=\exp(\int_0^1\log|P(e^{2\pi i t})|dt),\end{equation} 
the geometric mean of $|P(z)|$ for $z$ on the unit circle. However, $M(P)$ had appeared earlier in a paper by Lehmer \cite{Lehm33}, in an alternative form
\begin{equation} M(P)=|a_d|\prod\limits_{|\alpha_k|>1}|\alpha_k|.\end{equation}
The equivalence of the two definitions follows immediately from Jensen's formula \cite{Jen99}
\begin{equation*} \int_0^1\log|e^{2\pi i t}-\alpha|dt=\log_+|\alpha|,\end{equation*}
where $\log_+x$ denotes $\max(0,\,\log x).$ We will also deal with the \textit{small Mahler measure} which is defined as 
$$m(P):=\log M(P)=\int_0^1\log|P(e^{2\pi i t})|dt.$$
The concept of Mahler measure can be naturally extended to the class of Laurent polynomials $P(z)=a_{0}z^{p}+a_{1}z^{p+1}+\ldots+a_{d-1}z^{p+d-1}
+a_{d}z^{p+d}= a_{d}z^p\prod\limits_{k=1}^{d}(z-\alpha_k),$ where $a_0, a_d\ne0$ and $p$ is an arbitrary integer (not necessarily positive).
\section{Associated polynomials and their   properties}\label{polynomials}

The aim of this section is to introduce a few polynomials naturally associated with the circulant graph 
$$G_n=C_{\beta n}(s_1, \ldots,s_k,\alpha_1n,\ldots,\alpha_\ell n),\,
 1\le s_1<\ldots<s_k\le[\frac{\beta n}{2}],\, 1\le \alpha_1< \ldots<\alpha_\ell\le[\frac{\beta}{2}].$$ We start with the Laurent polynomial
$L(z)=2(k+l)-\sum\limits_{i=1}^k(z^{s_i}+z^{-s_i})-
\sum\limits_{m=1}^{\ell}(z^{\alpha_mn}+z^{-\alpha_mn})$  responsible for the structure of Laplacian of graph $G_n.$ 
More precisely, the Laplacian of $G_n$ is given by the matrix
$$\mathbb{L}=L(T)=2(k+l)I_{\beta n}-\sum\limits_{i=1}^k(T^{s_i}+T^{-s_i})
-\sum\limits_{m=1}^{\ell}(T^{\alpha_mn}+T^{-\alpha_mn}), $$   
 where $T$ the circulant matrix 
$circ(\underbrace{0,1,\ldots,0}_{\beta n}).$  We decompose $L(z)$ into the sum of two polynomials
$L(z)=P(z)+p(z^n),$  where $P(z)=2k-\sum\limits_{i=1}^k(z^{s_i}+z^{-s_i})$  and $p(z)=2l-
\sum\limits_{m=1}^{\ell}(z^{\alpha_m}+z^{-\alpha_m}).$    Now, we have to introduce a family of Laurent polynomials differing by a constant from $P(z).$  They are $P_u(z)=P(z)+p(e^{\frac{2\pi\,i\,u}{\beta}}),\,u=0,1,\ldots,\beta-1.$ One can check  that $P_{u}(z)=2k-\sum\limits_{i=1}^k(z^{s_i}+z^{-s_i})+4\sum\limits_{m=1}^{\ell}\sin^2(\frac{\pi \,u\,\alpha_{m}}{\beta}).$ In particular, $P_{0}(z)=P(z).$

We note that all the above polynomials are {\it palindromic}, that is they are invariant under replacement  $z$ by $1/z.$  A  non-trivial   palindromic Laurent polynomial can be represented in the form $\mathcal{P}(z)=a_s z^{-s}+a_{s-1}z^{-(s-1)}+\ldots+a_0+\ldots+a_{s-1}z^{s-1}+a_s z^s,$  where $a_s\neq 0.$  We will refer to $2s$  as a {\it degree}  of the polynomial $\mathcal{P}(z).$  Since $\mathcal{P}(z)=\mathcal{P}(\frac{1}{z}),$  the following polynomial of degree $s$  is well defined
$$\mathcal{Q}(w)=\mathcal{P}(w+\sqrt{w^2-1}).$$  We will call it a {\it Chebyshev trasform} of  $\mathcal{P}(z).$  Since $T_k(w)=\frac{(w+\sqrt{w^2-1})^k+(w+\sqrt{w^2-1})^{-k}}{2}$  is the Chebyshev polynomial of the first kind, one can easy deduce that
$$\mathcal{Q}(w)=a_0+2a_{1}T_1(w)+\ldots+2a_{s-1}T_{s-1}(w)+2a_{s}T_s(w).$$  Also, we  have  $\mathcal{P}(z)=\mathcal{Q}(\frac12(z+\frac{1}{z})).$  

Throughout the paper, we will use the following observation.   If $z_1,1/z_1,\ldots,z_s,1/z_s$  is the full list of the roots of $\mathcal{P}(z),$ then $w_k=\frac12(z_k+\frac{1}{z_k}),\,k=1,2,\ldots, s$  are all roots of the polynomial $\mathcal{Q}(w).$

By direct calculation, we obtain that  the Chebyshev transform of   polynomial $P_u(z)$  is 
$$Q_u(w)=2k-2\sum\limits_{i=1}^kT_{s_i}(w)+4\sum\limits_{m=1}^{\ell}\sin^2(\frac{\pi \,u\,\alpha_{m}}{\beta}).$$

In particular, if $z_{s}(u),\,1/z_{s}(u),\,s=1,2,\ldots,s_k$ are the roots of $P_{u}(z),$ then $w_{s}(u)=\frac12(z_{s}(u)+z_{s}(u)^{-1}),\,s=1,2,\ldots,s_k$ are all roots of the algebraic equation
$\sum_{i=1}^{k}T_{s_{i}}(w)=k+2\sum_{m=1}^{\ell}\sin^2(\frac{\pi\,u\,\alpha_{m}}{\beta}).$ To find asymptotic behavior for the number of spanning trees in the graph $G_n$ we also need the following lemma.
\bigskip

\begin{lemma}\label{lemma1}  Let $
\gcd(\alpha_1,\alpha_2,\ldots,\alpha_{\ell},\beta)=1.$ Suppose that  $P_u(z)=0,$  where $0<u<\beta.$
Then $|z|\neq 1.$
\end{lemma}
\begin{pf} 
First of all, we show that $p(e^{\frac{2\pi\,i\,u}{\beta}})=4\sum\limits_{m=1}^{\ell}\sin^2(\frac{\pi \,u\,\alpha_{m}}{\beta})>0.$ Indeed, suppose that $p(e^{\frac{2\pi\,i\,u}{\beta}})=0.$ Then there are integers $m_{j}$ such that $u\,\alpha_{j}=m_{j}\beta,\,j=1,2,\ldots,\ell.$
Hence 
$$B=\gcd(u\,\alpha_{1},\ldots,u\,\alpha_{\ell},u\,\beta)=u\gcd(\alpha_{1},\ldots,\alpha_{\ell},\beta)=u<\beta.$$From the other side  
$$B=\gcd(m_{1}\beta,\ldots,m_{\ell}\beta,u\,\beta)=\beta\gcd(m_{1},\ldots,m_{\ell},u)\geq\beta.$$
Contradiction.
Now, let $|z|=1.$  Then $z=e^{i \varphi},\text{ for some  } \varphi\in\mathbb{R}.$  We have $P_u(z)=P(z)+p(e^{\frac{2\pi\,i\,u}{\beta}}),$
where $P(z)=2k-\sum\limits_{i=1}^k(z^{s_i}+z^{-s_i})=2\sum\limits_{i=1}^{k}(1-\cos(s_{i}\varphi))\geq0.$ Hence, $P_u(z)>0$ and lemma is proved.
\end{pf}
\section{Complexity of circulant graphs with non-fixed jumps}\label{count}

The aim of this section is to find  new formulas for the numbers of spanning trees of circulant graph $C_{\beta n}(s_1,s_2,\ldots,s_k,\alpha_1n,\alpha_2n,\ldots,\alpha_\ell n)$  in terms of Chebyshev polynomials. It should be noted that nearby results were obtained earlier by different methods in the papers \cite{GolinXY2010, Louis, XiebinLinZhang, ZhangYongGol, ZhangYongGolin, MZXX2015}.

\begin{theorem}\label{theorem1}
The number of spanning trees in the circulant graph with non-fixed jumps
$$C_{\beta n}(s_1,\ldots,s_k,\alpha_1n,\ldots,\alpha_\ell n), 
\,1\le s_1<\ldots<s_k\leq[\frac{\beta n}{2}],\, 1\le \alpha_1<
\ldots <\alpha_\ell\leq[\frac{\beta}{2}]$$ is given by the formula
$$\tau(n)=\frac{n}{\beta\,q}\prod_{u=0}^{\beta-1}\prod_{j=1,\atop w_{j}(0)\neq1}^{s_k}
|2T_n(w_{j}(u))-2\cos(\frac{2\pi u}{\beta})|,$$
where for each $u=0,1,\ldots,\beta-1$ the numbers $w_{j}(u),\,j=1,2,\ldots,s_{k},$ are all the roots 
of the equation $\sum_{i=1}^k T_{s_i}(w)=k+2\sum_{m=1}^{\ell}\sin^2(\frac{\pi\,u\,\alpha_{m}}{\beta})
,\,T_{s}(w)$ is the Chebyshev polynomial of the first kind  and $q=s_1^2+s_2^2+\ldots+s_k^2.$
\end{theorem}

\medskip\noindent{Proof:} Let $G=C_{\beta n}(s_1,s_2,\ldots,s_k,\alpha_1n,\alpha_2n,
\ldots,\alpha_\ell n).$ By the celebrated Kirchhoff theorem, the number of spanning 
trees $\tau(n)$ in $G_n$ is equal to the product of non-zero eigenvalues of the 
Laplacian of a graph $G_n$ divided by the number of its vertices $\beta n.$ To 
investigate the spectrum of Laplacian matrix, we denote by $T$ the circulant matrix 
$circ(\underbrace{0,1,\ldots,0}_{\beta n}).$ Consider the Laurent polynomial 
$L(z)=2(k+l)-\sum\limits_{i=1}^k(z^{s_i}+z^{-s_i})-
\sum\limits_{m=1}^{\ell}(z^{\alpha_mn}+z^{-\alpha_mn}).$
Then the Laplacian of $G_n$ is given by the matrix
$$\mathbb{L}=L(T)=2(k+l)I_{\beta n}-\sum\limits_{i=1}^k(T^{s_i}+T^{-s_i})
-\sum\limits_{m=1}^{\ell}(T^{\alpha_mn}+T^{-\alpha_mn}).$$   

The eigenvalues of the circulant matrix $T$ are $\zeta_{\beta n}^j,\,j=0,1,\ldots,\beta n-1,$ 
where $\zeta_{m}=e^\frac{2\pi i}{m}.$ Since all of them are distinct, the matrix $T$ 
is conjugate to the diagonal matrix 
$\mathbb{T}=diag(1,\zeta_{\beta n},\ldots,\zeta_{\beta n}^{\beta n-1}).$ To find spectrum of $\mathbb{L},$ without loss of generality, one can assume that $T=\mathbb{T}.$
Then $\mathbb{L}$ is a diagonal matrix. This essentially simplifies the problem of finding 
eigenvalues of $\mathbb{L}.$ Indeed, let $\lambda$ be an eigenvalue of $\mathbb{L}$ and
$x$ be the respective eigenvector. Then we have the following system of linear equations
$$((2(k+l)-\lambda)I_{\beta n}-\sum\limits_{i=1}^k({T}^{s_i}+{T}^{-s_i})
-\sum\limits_{m=1}^{\ell}({T}^{\alpha_{m}n}+{T}^{-\alpha_{m}n}))x=0.$$

Let $\textbf{e}_{j} =(0,\ldots,\underbrace{1}_{j-th},\ldots, 0),\,j=1,\ldots,\beta n.$
The $(j,j)$-th entry of $\mathbb{T}$ is equal to $\zeta_{\beta n}^{j-1}.$
Then, for $j=0,\ldots, \beta n-1,$ the matrix $\mathbb{L}$ has an eigenvalue
\begin{equation}\label{lambda}\lambda_j=L(\zeta_{\beta n}^j)=
2(k+l)-\sum\limits_{i=1}^k(\zeta_{\beta n}^{j s_i}+\zeta_{\beta n}^{-j s_i})
-\sum\limits_{m=1}^{\ell}(\zeta_{\beta}^{j\alpha_m}+\zeta_{\beta}^{-j\alpha_m}),\end{equation}
with eigenvector $\textbf{e}_{j+1}.$ Since all graphs under consideration are supposed to be connected, we have $\lambda_0=0$ and 
$\lambda_j>0,\, j=1,2,\ldots,\beta n-1.$ Hence
\begin{equation}\label{kirchhoff}\tau(n)=\frac{1}{\beta n}
\prod\limits_{j=1}^{\beta n-1}L(\zeta_{\beta n}^j).\end{equation}

By setting $j=\beta t+u,$ where $0\le t\le n-1,\,0\le u\le\beta-1,$
we rewrite the formula (\ref{kirchhoff}) in the form

\begin{equation}\label{product}\tau(n)=(\frac{1}{n}\prod\limits_{t=1}^{n-1}
L(\zeta_{\beta n}^{\beta t}))(\frac{1}{\beta}
\prod\limits_{u=1}^{\beta-1}\prod\limits_{t=0}^{n-1}L(\zeta_{\beta n}^{ t \beta + u})).\end{equation}

It is easy to see that $\tau(n)$ is the product of two numbers $\tau_{1}(n)=\frac{1}{n}\prod\limits_{t=1}^{n-1}L(\zeta_{\beta n}^{\beta t})$ and $\tau_{2}(n)=\frac{1}{\beta}
\prod\limits_{u=1}^{\beta-1}\prod\limits_{t=0}^{n-1}L(\zeta_{\beta n}^{ t \beta + u}).$

We note that $$L(\zeta_{\beta n}^{\beta t}) = 2k - \sum\limits_{i=1}^k(\zeta_{\beta n}^{\beta t s_i} +\zeta_{\beta n}^{-\beta t s_i}) = 2k - \sum\limits_{i=1}^k(\zeta_{n}^{t s_i} + \zeta_{n}^{-t s_i})=P(\zeta_{n}^{t}),\,1\le t \le n-1.$$ The numbers $\mu_t=P(\zeta_{n}^{t}),\,1\le t\le n-1$ run through all non-zero eigenvalues of circulant graph $C_{n}(s_1,s_2,\ldots,s_k)$ with fixed jumps $s_1,s_2,\ldots,s_k$ and $n$ vertices. So $\tau_{1}(n)$ coincide with the number of spanning trees in $C_{n}(s_1,s_2,\ldots,s_k).$ By (\cite{MedMed3}, Corollary~1) we get 
\begin{equation}\label{tau1}\tau_{1}(n)=\frac{n}{q}\prod_{j=1,\atop w_{j}(0)\neq1}^{s_k}|2T_n(w_{j}(0))-2|,\end{equation} where  $w_{j}(0),\,j=1,2,\ldots,s_k,$ are all the roots of the equation $\sum_{i=1}^k T_{s_i}(w)=k.$ 

In order to continue the calculation of $\tau(n)$ we have to find the product 
$$\tau_{2}(n)=\frac{1}{\beta}
\prod\limits_{u=1}^{\beta-1}\prod\limits_{t=0}^{n-1}L(\zeta_{\beta n}^{ t \beta + u}).$$ 
Recall that $L(z)=P(z)+p(z^n).$  
Since 
$(\zeta_{\beta n}^{ \beta t + u})^n=\zeta_{\beta}^{ \beta t + u}=\zeta_{\beta}^{u},$ we obtain 

$$L(\zeta_{\beta n}^{\beta t+ u})=P(\zeta_{\beta n}^{\beta t + u})+p(\zeta_{\beta }^{\beta t+ u})=P(\zeta_{\beta n}^{\beta t+ u})+p(\zeta_{\beta }^{u})=P_{u}(\zeta_{\beta n}^{\beta t + u}),$$ where
$P_{u}(z)=P(z)+p(\zeta_{\beta }^{u}).$  By Section~\ref{polynomials},  we already know that

$P_u(z)=(-1)\prod\limits_{j=1}^{s_k}(z-z_{j}(u))(z-z_{j}(u)^{-1}),$ where $w_{j}(u)=\frac12(z_{j}(u)+z_{j}(u)^{-1}),\,j=1,2,\ldots,s_k$ are all roots of the equation
$\sum_{i=1}^k T_{s_i}(w)=k+2\sum_{d=1}^{\ell}\sin^2(\frac{\pi\,u\,\alpha_{d}}{\beta}).$

We note that $\zeta_{\beta n}^{ t \beta + u}=e^{\frac{i(2\pi t+\omega_u)}{n}},$  where  
$\omega_u=\frac{2\pi u}{\beta}.$      
Then $\prod\limits_{t=0}^{n-1}L(\zeta_{\beta n}^{ t \beta + u})=\prod\limits_{t=0}^{n-1}P_u(e^{\frac{i(2\pi t+\omega_u)}{n}}).$      


To evaluate the latter product, we need  following lemma.

\begin{lemma}\label{lemmaH} Let $H(z)=\prod\limits_{s=1}^{m}(z-z_s)(z-z_s^{-1})$ and 
$\omega$ is a real number. Then
$$\prod\limits_{t=0}^{n-1}H(e^{\frac{i(2\pi t+\omega)}{n}})=(-e^{i\omega})^{m}
\prod\limits_{s=1}^{m}(2T_n(w_s)-2\cos(\omega)),$$ 
where $w_s=\frac12(z_s+z_s^{-1}),\,s=1,\ldots,k$ and $T_n(w)$ is the Chebyshev 
polynomial of the first kind.
\end{lemma}

\begin{pf} We note that $\frac{1}{2}(z^{n}+z^{-n})=T_n(\frac{1}{2}(z+z^{-1})).$ 
By the substitution $z=e^{i\,\varphi},$ this follows from the evident identity 
$\cos(n\varphi)=T_n(\cos\varphi).$ Then we have
\begin{eqnarray*}
\prod\limits_{t=0}^{n-1}H(e^{\frac{i(2\pi t+\omega)}{n}})&=& 
\prod\limits_{t=0}^{n-1}\prod\limits_{s=1}^{m}(e^{\frac{i(2\pi t+\omega)}{n}}-z_s)
(e^{\frac{i(2\pi t+\omega)}{n}}-z_s^{-1})\\
&=&\prod\limits_{s=1}^{m}\prod\limits_{t=0}^{n-1}(-e^{\frac{i(2\pi t+\omega)}{n}}z_{s}^{-1})
(z_{s}-e^{\frac{i(2\pi t+\omega)}{n}})(z_{s}-e^{-\frac{i(2\pi t+\omega)}{n}})\\
&=&\prod\limits_{s=1}^{m}(-e^{i\omega}{z_{s}}^{-n})\prod\limits_{t=0}^{n-1}
(z_{s}-e^{\frac{i(2\pi t+\omega)}{n}})(z_{s}-e^{-\frac{i(2\pi t+\omega)}{n}})\\
&=&\prod\limits_{s=1}^{m}(-e^{i\omega}{z_{s}}^{-n})(z_{s}^{2n}-2\cos(\omega)z_{s}^{n}+1)\\
&=&\prod\limits_{s=1}^{m}(-e^{i\omega})(2\,\frac{z_{s}^{n}+z_{s}^{-n}}{2}-2\cos(\omega))\\
&=&(-e^{i\omega})^{m}\prod\limits_{s=1}^{m}(2T_n(w_s)-2\cos(\omega)).
\end{eqnarray*}
\end{pf}
\bigskip

Since $P_u(z)=-H_u(z),$ where $H_u(z)=\prod\limits_{j=1}^{s_k}(z-z_{j}(u))(z-z_{j}(u)^{-1}),$ by Lemma~\ref{lemmaH} we get

$$\prod\limits_{t=0}^{n-1}P_{u}(e^{\frac{i(2\pi t+\omega_u)}{n}})=(-1)^{n}(-e^{\frac{2\pi\,i\,u}{\beta}})^{s_{k}}\prod\limits_{j=1}^{s_{k}}(2T_{n}(w_{j}(u))-2\cos(\frac{2\pi\,u}{\beta})).$$    Then,

\begin{eqnarray}\label{tausign}
\nonumber\tau_{2}(n)&=&\frac{1}{\beta}\prod\limits_{u=1}^{\beta-1}
\prod\limits_{t=0}^{n-1}L(\zeta_{\beta n}^{\beta t+u})=\frac{1}{\beta}\prod\limits_{u=1}^{\beta-1}
\prod\limits_{t=0}^{n-1}P_{u}(e^{\frac{i(2\pi j+\omega_u)}{n}})\\
&=&\frac{(-1)^{n(\beta-1)}}{\beta}\prod\limits_{u=1}^{\beta-1}(-e^{\frac{2\pi\,i\,u}{\beta}})^{s_{k}}
\prod\limits_{j=1}^{s_{k}}(2T_{n}(w_{j}(u))-2\cos(\frac{2\pi\,u}{\beta}))\\
\nonumber&=&\frac{(-1)^{n(\beta-1)}}{\beta}\prod\limits_{u=1}^{\beta-1}
\prod\limits_{j=1}^{s_{k}}(2T_{n}(w_{j}(u))-2\cos(\frac{2\pi\,u}{\beta})).
\end{eqnarray}
Since the number $\tau_2(n)$  is a product of positive eigenvalues of $G_n$ divided by $\beta$, from (\ref{tausign}) we have
\begin{equation}\label{tau2}
\tau_2(n)=\frac{1}{\beta}\prod\limits_{u=1}^{\beta-1}
\prod\limits_{j=1}^{s_{k}}|2T_{n}(w_{j}(u))-2\cos(\frac{2\pi\,u}{\beta})|.
\end{equation}

Combining   equations (\ref{tau1})  and (\ref{tau2})  we finish the proof of the theorem.
$\hfill \qed$

\bigskip

As the first consequence from   Theorem~\ref{theorem1} we have the following result  obtained earlier by Justine Louis \cite{Louis}  in a slightly different form.

\begin{corollary}\label{louis1} The number of spanning trees in the  circulant graphs with non-fixed jumps
$C_{\beta n}(1,\alpha_1n,\alpha_2n,\ldots,\alpha_\ell n),$ where
$ 1\le \alpha_1<\alpha_2<\ldots<\alpha_\ell\leq[\frac{\beta}{2}]$ is given by the formula
$$\tau(n)=\frac{n\,2^{\beta-1}}{\beta}\prod_{u=1}^{\beta-1}(T_n(1+2\sum_{m=1}^{\ell}\sin^2(\frac{\pi\,u\,\alpha_{m}}{\beta}))-\cos(\frac{2\pi\,u}{\beta})),$$ where $T_n(w)$ is the Chebyshev polynomial of the first kind.
\end{corollary}

\begin{pf} Follows directly from the theorem. \end{pf}

The next corollary is new.

\begin{corollary}\label{louis2} The number of spanning trees in the  circulant graphs with non-fixed jumps 
$C_{\beta n}(1,2,\alpha_1n,\alpha_2n,\ldots,\alpha_\ell n),$ where 
$ 1\le \alpha_1<\alpha_2<\ldots<\alpha_\ell\leq[\frac{\beta}{2}]$ is given by the formula 
$$\tau(n)=\frac{n F_n^2}{\beta}\,\prod_{u=1}^{\beta-1}
\prod_{j=1}^2|2T_n(w_{j}(u))-2\cos(\frac{2\pi\,u}{\beta})|,$$ where $F_n$  
is the $n$-th Fibonacci number, $T_n(w)$ is the Chebyshev polynomial of the first kind and $w_{1,2}(u) =\left(-1\pm\sqrt{25+16\sum_{m=1}^{\ell}\sin^2(\frac{\pi\,u\,\alpha_{m}}{\beta})}\,\right)/4.$
   
\end{corollary}  
We note that $n F_n^2$ is the number of spanning trees in the graph $C_{n}(1,2).$

\begin{pf} 
In this case, $k=2,\,s_1=1,\,s_2=2$ and $q=s_1^2+s_2^2=5.$ Given $u$ we   find $w_{j}(u),\,j=1,2$ as the roots of the algebraic equation 
$$T_{1}(w)+T_2(w)=2+2\sum_{m=1}^{\ell}\sin^2(\frac{\pi\,u\,\alpha_{m}}{\beta}),$$ 
where $T_1(w)=w$ and $T_2(w)=2w^2-1.$ For $u=0$ the roots are $w_1(0)=1$  and $w_2(0)=-3/2.$ Hence $\tau_1(n)=\frac{n}{5}|2T_n(-\frac{3}{2})-2|=\frac{n}{5}|(\frac{-3+\sqrt{5}}{2})^n+(\frac{-3-\sqrt{5}}{2})^n-2|=n F_n^2$ gives the well-known formula for the number of spanning trees in the graph $C_{n}(1,2).$ (See, for example, \cite{BoePro},\,Theorem~4). For $u>0$ the numbers $w_{1,2}(u)$ are   roots of the quadratic equation 
$$2w^2+w-3-2\sum_{m=1}^{\ell}\sin^2(\frac{\pi\,u\,\alpha_{m}}{\beta})=0.$$ 
By (\ref{tau2}) we get $\tau_2(n)=\frac{1}{\beta}\,\prod_{u=1}^{\beta-1}\prod_{j=1}^2|2T_n(w_{j}(u))-2\cos(\frac{2\pi\,u}{\beta})|.$ Since $\tau(n)=\tau_1(n)\tau_2(n),$ the result follows.
\end{pf}

\section{Arithmetic properties of the complexity for circulant graphs}\label{circarithm}

It was noted in the series of paper (\cite{KwonMedMed}, \cite{Med1}, \cite{MedMed2}, \cite{MedMed3}) that in many important cases the complexity of graphs is given by the formula $\tau(n)=p\, n\, a(n)^2,$ where $a(n)$ is an integer sequence and $p$ is a prescribed constant depending only of parity of $n.$   

The aim of the next theorem is to explain this phenomena for circulant graphs with non-fixed jumps. Recall that any positive integer $p$ can be uniquely represented in the form $p=q \,r^2,$ where $p$ and $q$ are positive integers and $q$ is square-free. We will call $q$ the \textit{square-free part} of $p.$
\begin{theorem}\label{theorem2} Let $\tau(n)$ be the number of spanning trees of the circulant graph
$$G_n=C_{\beta n}(s_1,s_2,\ldots,s_k,\alpha_1n,\alpha_2n,\ldots,\alpha_\ell n),$$
where $1\le s_1<s_2<\ldots<s_k\leq[\frac{\beta n}{2}],
\,1\le\alpha_1<\alpha_2<\ldots,\alpha_\ell\le[\frac{\beta}{2}].$

Denote by $p$ and $q$ the number of odd elements in the sequences $s_1,s_2,\ldots,s_k$
and $\alpha_1,\alpha_2,\ldots,\alpha_\ell$ respectively. Let $r$ be the square-free part of $p$
and $s$ be the square-free part of $p+q.$ Then there exists an integer sequence $a(n)$ such that

\begin{enumerate}
\item[ $1^0$ ]  $\tau(n)=\beta\, n\,a(n)^2,$ if $n$ and  $\beta$ are odd;
\item[ $2^0$ ]  $\tau(n)=\beta\,r\,n\,a(n)^2,$  if $n$ is even;
\item[ $3^0$ ]  $\tau(n)=\beta\,s\,n\,a(n)^2,$  if  $n$ is odd and $\beta$ is even.
\end{enumerate}
\end{theorem}

\begin{pf} The number of odd elements in the sequences $s_1,s_2,\ldots,s_k$
and $\alpha_1,\alpha_2,\ldots,\alpha_\ell,$ respectively is counted respectively by the formulas
$p=\sum\limits_{i=1}^{k}\frac{1-(-1)^{s_i}}{2}$    and $q=\sum\limits_{i=1}^{\ell}\frac{1-(-1)^{\alpha_i}}{2}.$
We already know that all non-zero eigenvalues of the graph $G_n$
are given by the formulas
$\lambda_j=L(\zeta_{\beta n}^{j}),\,j=1,\ldots,\beta n-1,$ where  $\zeta_{\beta n}=e^{\frac{2\pi i}{\beta n}}$ and 
\begin{equation*}\label{lambdaz} L(z)=
2(k+l)-\sum\limits_{i=1}^k(z^{s_i}+z^{-s_i})
-\sum\limits_{m=1}^{\ell}(z^{n\alpha_m}+z^{-n \alpha_m }).
\end{equation*} We note that $\lambda_{ \beta n-j}=L(\zeta_{\beta n}^{\beta n-j})=L(\zeta_{\beta n}^{j})=\lambda_j.$

By the Kirchhoff theorem we have $\beta n\,\tau(n)=\prod\limits_{j=1}^{\beta n-1}\lambda_j.$
Since $\lambda_{\beta n-j}=\lambda_j,$  we obtain
$\beta n\,\tau(n)=(\prod\limits_{j=1}^{\frac{\beta n-1}{2}}\lambda_j)^2$ if $\beta n$ is odd and
$\beta n\,\tau(n)=\lambda_{\frac{\beta n}{2}}(\prod\limits_{j=1}^{\frac{\beta n}{2}-1}\lambda_j)^2$
if $\beta n$ is even. We note that each algebraic number $\lambda_j$ comes with all its
Galois conjugate \cite{Lor}. So, the numbers
$c(n)=\prod\limits_{j=1}^{\frac{\beta n-1}{2}}\lambda_j$ and $d(n)=\prod\limits_{j=1}^{\frac{\beta n}{2}-1}\lambda_j$
are integers. Also, for even $n$ we have 
$$\lambda_{\frac{\beta n}{2}}=L(-1)=2(k+l)-\sum\limits_{i=1}^k((-1)^{s_i}+(-1)^{-s_i})
-\sum\limits_{m=1}^{\ell}((-1)^{n\alpha_m}+(-1)^{-n \alpha_m })$${}
$$=2k-\sum\limits_{i=1}^k((-1)^{s_i}+(-1)^{-s_i})=4\sum\limits_{i=1}^{k}\frac{1-(-1)^{s_i}}{2} =4p.$$

If $n$ is odd and $\beta$  is even,  the number $\frac{\beta n}{2}$  is integer again. Then we obtain
$$\lambda_{\frac{\beta n}{2}}=L(-1)=2(k+l)-\sum\limits_{i=1}^k((-1)^{s_i}+(-1)^{-s_i})
-\sum\limits_{m=1}^{\ell}((-1)^{\alpha_m}+(-1)^{-\alpha_m })$${}
$$=4\sum\limits_{i=1}^{k}\frac{1-(-1)^{s_i}}{2}+4\sum\limits_{m=1}^{\ell}\frac{1-(-1)^{\alpha_m}}{2}=4p+4q.$$
Therefore, $\beta\,n\,\tau(n)= c(n)^2$ if $\beta$ and $ n$ are odd, $\beta\,n\,\tau(n)=4 p\,d(n)^2$ if $n$ is even and $\beta\,n\,\tau(n)=4(p+q)\,d(n)^2$ if $n$ is odd and $\beta$ is even. Let $r$ be the square-free part of $p$ and $s$ be the square-free part of $p + q.$ Then there are integers $u$ and $v$ such that $p=ru^2$ and $s=(p+q) v^2.$ Hence,
\begin{enumerate}
\item[$1^{\circ}$] $\displaystyle{\frac{\tau(n)}{\beta\,n}=\left(\frac{c(n)}{\beta\,n}\right)^2}$  if $n$  and  $\beta$ are odd, 
\item[$2^{\circ}$]  $\displaystyle{\frac{\tau(n)}{\beta\,n}=r\left(\frac{2\,u\,d(n)}{\beta\,n}\right)^2}$ if $n$  is even and
\item[$3^{\circ}$]  $\displaystyle{\frac{\tau(n)}{\beta\,n}=s\left(\frac{2\,v\,d(n)}{\beta\,n}\right)^2}$    if $n$ is odd and $\beta$  is even.
\end{enumerate}

Consider an automorphism group $\mathbb{Z}_{\beta\,n}$ of the graph $G_n$   consisting of elements circularly permuting its vertices $v_0,v_1,\ldots,v_{\beta\,n-1}$ and acting without fixed edges. Such a group always exists,    since in the case of even $\beta\,n$  we have {\it even} number of multiple edges between the opposite vertices $v_i$  and  $v_{i+\frac{\beta\,n}{2}},$ where the indices are taken $\textrm{mod}\,{\beta\,n}.$ 

The group $\mathbb{Z}_{\beta\,n}$ acts  fixed point free on the set vertices of $G_n.$ We are aimed to show that it also acts    freely on the set of the spanning trees in   the graph.  Indeed, suppose that some non-trivial element $\gamma$  of $\mathbb{Z}_{\beta\,n}$  leaves a spanning tree $A$  in the graph $G_n$ invariant. Then $\gamma$ fixes the center  of  $A.$ The center of a tree is a vertex or an edge. The first case is impossible, since  $\gamma$ acts freely on the set of vertices.  In the second case, $\gamma$  permutes
the endpoints of an edge connecting the apposite 
vertices of $G_n.$ This means that   $\beta\,n$ is even, and $\gamma$ is the unique involution in the group $\mathbb{Z}_{\beta\,n}.$    This is also impossible, since  the group is acting without fixed edges.
 
So, the cyclic group  $\mathbb{Z}_{\beta n}$  acts  on the set of  spanning trees of the graph $G_{n}$ fixed point free. Therefore $\tau(n)$ is a
multiple of $\beta\,n$  and their quotient  $\frac{\tau(n)}{\beta\,n}$ is an integer.

Setting $a(n)=\frac{c(n)}{\beta\,n}$ in the  case $1^{\circ},\ a(n)=\frac{2\,u\,d(n)}{\beta\,n}$ in the case $2^{\circ},$ and $a(n)=\frac{2\,v\,d(n)}{\beta\,n}$ in the case $3^{\circ},$ we conclude that number $a(n)$ is always integer and the statement of the theorem follows.
\end{pf}
\bigskip

\section{Asymptotic for the number of spanning trees}\label{assection}

In this section, we give asymptotic formulas for the number of spanning trees for circulant graphs.
It is interesting to compare these results with those in papers \cite{GolinXY2010, Louis, XiebinLinZhang, ZhangYongGol, MZXX2015}, where the similar results were obtained by different methods.

\begin{theorem}\label{asymptotic1}
The number of spanning trees in the circulant graph
\begin{center}$C_{\beta n}(s_1,s_2,\ldots,s_k,\alpha_1n,\alpha_2n,\ldots,\alpha_\ell n),$\end{center}
 $1\le s_1<s_2<\ldots<s_k\leq~[\frac{\beta n}{2}],
\,1\le \alpha_1<\alpha_2<\ldots<\alpha_\ell\le[\frac \beta2],$\,\ ${\rm gcd} (s_1,s_2,\ldots,s_k)=d,\,
{\rm gcd}(\alpha_1,\alpha_2,\ldots,\alpha_\ell, \beta)=\delta,$ and $d$ and $\delta$  are relatively prime
has the following asymptotic
$$\tau(n)\sim\frac{ n\, \delta^2\,d^2}{\beta\,q}A^n,\text{ as }n\to\infty \text{  and   }(n,d)=1,$$
where $q=s_1^2+s_2^2+\ldots+s_k^2,\,A=\prod_{u=0}^{\beta-1} M(P_u)$ and $M(P_u)=\exp(\int_0^1\log|P_u(e^{2\pi i t})|dt)$
is the Mahler measure of Laurent polynomial $P_u(z)=2k-\sum\limits_{i=1}^k(z^{s_i}+z^{-s_i})+4\sum\limits_{m=1}^{\ell}\sin^2(\frac{\pi u \alpha_m }{\beta}).$
\end{theorem}\begin{pf} Without loss of generality, we can restrict ourself by the case $\delta = 1.$ Indeed, 
if $\delta > 1,$ then  one can replace the graph $ C_{\beta n}(s_1,s_2,\ldots,s_k,\alpha_1n,\alpha_2n,\ldots,\alpha_\ell n)$  by an isomorphic graph $$G_m^{\prime}=C_{\beta^{\prime} m}(s_1,s_2,\ldots,s_k,\alpha^{\prime}_1m,\alpha^{\prime}_2m,\ldots,\alpha^{\prime}_\ell m),$$  where
$\beta^{\prime}=\beta/\delta,\,\alpha^{\prime}_1=\alpha_1/\delta,\,\ldots,\alpha^{\prime}_{\ell}=\alpha_{\ell}/\delta,\text{ and }\,m=\delta\,n,$  with $\gcd(\alpha^{\prime}_1,\alpha^{\prime}_2,\ldots,\alpha^{\prime}_{\ell},\beta^{\prime})=1.$ 
From now on, we suppose that ${\rm gcd}(\alpha_1,\alpha_2,\ldots,\alpha_\ell, \beta)=1.$

By Theorem~\ref{theorem1},  $\tau(n)=\tau_1(n)\tau_2(n),$ where $\tau_1(n)$ is the number of spanning trees in   $C_{n}(s_1,s_2,\ldots,s_k)$ and  $ \tau_2(n)=\frac{1}{\beta}\prod\limits_{u=1}^{\beta-1}
\prod\limits_{j=1}^{s_{k}}|2T_{n}(w_{j}(u))-2\cos(\frac{2\pi\,u}{\beta})|.$  By (\cite{MedMed3}, Theorem~5) we already  know that
$$\tau_1(n)\sim\frac{ n\, d^2}{q}A_0^n,\text{ as }n\to\infty \text{  and   }(n,d)=1,$$
where $A_0$ is the Mahler measure of Laurent polynomial $P_0(z).$ So, we have to find asymptotics for 
$\tau_2(n)$ only.

By Lemma~\ref{lemma1}, for any integer $u,\,0<u<\beta$ we obtain $T_n(w_j(u))=\frac{1}{2}(z_j(u)^n+z_j(u)^{-n}),$ where the $z_j(u)$ and $1/z_j(u)$ are roots of the polynomial $P_u(z)$ satisfying the inequality $|z_j(u)|\neq1,\,j=1,2,\ldots,s_{k}.$ Replacing $z_j(u)$ by $1/z_j(u),$ if necessary, we can assume that  $|z_j(u)|>1$ for all $j=1,2,\ldots,s_{k}.$ Then $T_{n}(w_j(u))\sim\frac{1}{2}z_{j}(u)^n,$ as $n$ tends to $\infty.$ So $|2T_n(w_j(u))-2\cos(\frac{2\pi\,u}{\beta})|\sim |z_j(u)|^n,\,\,n\to\infty.$ 
Hence
$$\prod_{j=1}^{s_k}|2T_n(w_j(u))-2\cos(\frac{2\pi\,u}{\beta})|\sim \prod_{s=1}^{s_k}|z_j(u)|^n=
 \prod\limits_{P_u(z)=0,\,|z|>1}|z|^n= A_u^n,$$  where
$A_u=\prod\limits_{P_u(z)=0,\,|z|>1}|z|$  coincides with the Mahler   measure of $P_u(z).$ As a result,
$$ \tau_2(n)=\frac{1}{\beta}\prod\limits_{u=1}^{\beta-1}
\prod\limits_{j=1}^{s_{k}}|2T_{n}(w_{j}(u))-2\cos(\frac{2\pi\,u}{\beta})|\sim\frac{1}{\beta}\prod\limits_{u=1}^{\beta-1}A_u^n.$$  

Finally, $ \tau(n)=\tau_1(n)\tau_2(n)\sim\frac{ n\, d^2 }{\beta\,q}\prod\limits_{u=0}^{\beta-1}A_u^n,\text{ as }n\to\infty \text{  and   }(n,d)=1.$
Since, $A_u=M(P_u),$  the result follows.

\end{pf}
\bigskip

As an immediate consequence of above theorem we have the following result obtained earlier in (\cite{GolinXY2010},\,Theorem~3) by completely different methods.
\bigskip

\begin{corollary}\label{mahler}
The thermodynamic limit of the sequence $C_{\beta n}(s_1,s_2,\ldots,s_k,\alpha_1n,\alpha_2n,\ldots,\alpha_\ell n)$
of circulant graphs is equal to the arithmetic mean of small Mahler measures of Laurent polynomials $P_u(z),\,u=0,1,\ldots,\beta-1.$ More precisely,
$$\lim\limits_{n\to\infty}\frac{\log\tau(C_{\beta n}(s_1,s_2,\ldots,s_k,\alpha_1n,\alpha_2n,\ldots,\alpha_\ell n))}{\beta\,n}=\frac{1}{\beta}\sum_{u=0}^{\beta -1}m(P_u),$$
where $m(P_u)=\int\limits_0^1\log|P_u(e^{2\pi i t})|dt$  and  $P_u(z)=2k-\sum\limits_{i=1}^k(z^{s_i}+z^{-s_i})+4\sum\limits_{m=1}^{\ell}\sin^2(\frac{\pi\,u\, \alpha_m}{\beta}).$
\end{corollary}
\section{Examples}\label{tables}

\begin{enumerate}

\item{\textbf{Graph $C_{2n}(1,n).$}} (M\"obius ladder with double steps). By Theorem~\ref{theorem1}, we have
$\tau(n)=\tau(C_{2n}(1,n))=n\,(T_n(3)+1).$ Compare this result with (\cite{ZhangGolin2002}, Theorem~4).
Recall \cite{BoePro} that the number of spanninig trees in the M\"obius ladder with single steps is given by the formula $n\,(T_n(2)+1).$
\
\item{\textbf{Graph $C_{2n}(1,2,n).$}} We have
$\tau(n)=2n F_n^2|T_n(\frac{-1-\sqrt{41}}4)-1||T_n(\frac{-1+\sqrt{41}}4)-1|.$
By Theorem~\ref{theorem2}, one can find an integer sequence $a(n)$ such that $\tau(n)=2na(n)^2$ if $n$ is even
and $\tau(n)=na(n)^2$ if $n$ is odd.

\item{\textbf{Graph $C_{2n}(1,2,3,n).$}} Here
$\tau(n)=\frac{8 n}{7}(T_n(\theta_1)-1)(T_n(\theta_2)-1)\prod_{p=1}^3(T_n(\omega_p)+1),$
where $\theta_{1,2}=\frac{-3\pm\sqrt{-7}}4$ and $\omega_p,\,p=1, 2, 3$  are roots of the
cubic equation $2w^3+w^2-w-3=0.$ We have $\tau(n)=6 n a(n)^2$ is $n$ is odd and
$\tau(n)=4 n a(n)^2$ is $n$ is even. Also, $\tau(n)\sim\frac{n}{28}A^{n},\,n\to\infty,$
where $A\approx 42.4038.$

\item{\textbf{Graph $C_{3n}(1,n).$}} We have
$$\tau(n)=\frac{n}{3}(2\,T_n(\frac52)+1)^2=\frac{n}{3}((\frac{5+\sqrt{21}}{2})^n+(\frac{5-\sqrt{21}}{2})^n+1)^2.$$ See also (\cite{ZhangGolin2002}, Theorem~5).
We note that $\tau(n)=3n a(n)^2,$ where $a(n)$ satisfies the recursive relation $a(n)=6a(n-1)-6a(n-2)+a(n-3)$ with initial data $a(1)=2,\,a(2)=8,a(3)=37.$
 
\item{\textbf{Graph $C_{3n}(1,2,n).$}}  By Theorem~\ref{theorem1},  we obtain
$$\tau(n)=\frac{n}{3}
F_n^2(2\,T_n(\omega_1)+1)^2(2\,T_n(\omega_2)+1)^2,$$
where $\omega_{1,2}=\frac{-1\pm\sqrt{37}}{4}.$
By Theorem~\ref{theorem2},  $\tau(n)=3na(n)^2$ for some integer sequence $a(n).$

\item{\textbf{Graph $C_{6n}(1,n,3n).$}}  Now, we get
$$\tau(n)=\frac{n}{3}(2\,T_n(\frac{5}2)+1)^2(2\,T_n(\frac{7}2)-1)^2(T_n(5)+1).$$
For a suitable integer sequence $a(n),$  one has
 $\tau(n)=6n a(n)^2$ if $n$ is even and $\tau(n)=18 n a(n)^2$ if $n$ is odd. 
 
\item{\textbf{Graph $C_{12n}(1,3n,4n).$}}  In this case
$$\tau(n)=\frac{2 n}{3}T_n(2)^2(2\,T_n(\frac{5}2)+1)^2(T_n(3)+1)(4\,T_n(\frac{7}2)^2-3)^2(2\,T_n(\frac{9}2)-1)^2.$$ 
By Theorem~\ref{theorem2}, one can conclude that $\tau(n)=3 n a(n)^2$ if $n$ is even and $\tau(n)=6 n a(n)^2$ if $n$ is odd, for some  sequence $a(n)$ of even numbers.
\end{enumerate}

\section*{ACKNOWLEDGMENTS}

The results of this work were partially supported by the Russian Foundation for Basic Research, (grants 18-01-00420 and 18-501-51021). The results given in Sections 5 and 6 were supported by the Laboratory of Topology and Dynamics, Novosibirsk State Uni- versity (contract no. 14.Y26.31.0025 with the Ministry of Education and Science of the Russian Federation).

\end{document}